\documentclass[russcorr,a4paper,12pt]{article}

\usepackage[russian]{babel}

\usepackage{graphicx}
\usepackage{amsmath}
\usepackage{amsfonts}
\usepackage{amssymb}
\newtheorem{theorem}{Теорема}
\newtheorem{corollary}{Следствие}
\newtheorem{definition}{Определение}
\newtheorem{example}{Пример}

\newtheorem{proposition}{Предложение}
\newtheorem{remark}{Замечание}
\newenvironment{proof}[1][Доказательство]{\textbf{#1.} }{\ \rule{0.5em}{0.5em}}
\title{Разрешимые алгебры Ли, заданные произведением образующих, и некоторые их приложения}
\author{{Рутвиг Кампоамор-Штурсберг}\\
{Universidad Complutense de Madrid}\\
email: rutwig@ucmail.ucm.es}
\date{}

\begin{document}

\maketitle

\begin{center}
{\cyr Аннотация\newline}
\end{center}

{\cyr В настоящей работе мы обобщаем определение произведения
образующими на класс разрешимых алгебр Ли. Мы исследуем число
инвариантных функций коприсоединенного представления этих алгебр
посредством уравнений Маурера-Картана и приводим некоторые
приложения к структурам произведения на алгебрах.}

{\cyr

\section{Введение}

Разрешимые алгебры Ли составляют важный класс алгебр не только в
связи с многочисленными приложениями к геометрии и механике, но
также и для классификации полупрямых сумм алгебр Ли (как известно,
полупростой случай классифицирован). В общем случае классификация
не была проведена из-за отсутствия общих структур для разрешимых
алгебр Ли (таких, как форма Киллинга). В настоящей работе мы
исследуем разрешимые алгебры Ли специального вида, которые мы
называем произведениями образующими. Будучи изначально
определенной для нильпотентных алгебр, операция произведения
образующими может быть определена и в случае разрешимых алгебр Ли.
Это приводит к методу построения неразложимых алгебр Ли в больших
размерностях.\newline Мы исследуем также два приложения:
инварианты коприсоединенного представления и структуры
произведения на группах Ли (мы изучаем алгебры Ли групп Ли). Мы
покажем, что произведения образующими открывают естественный
способ изучения структур произведения и паракомплексных структур в
произвольной размерности.

\section{Произведение образующими разрешимых алгебр Ли}

Произведение образующими было изначально введено для нильпотентных
алгебр для изучения систем весов. Введенные посредством графа
весов, оно может быть охарактеризовано расширениями, что позволяет
эффективно работать с уравнениями Маурера-Картана \cite{AC4}. Эта
операция может быть обобщена на случай разрешимых алгебр Ли; в
этом случае мы уже теряем наглядную интерпретацию посредством
графов. Эта операция позволяет строить неразложимые алгебры Ли.
\medskip

Как обычно, обозначим через $Z\left(  \frak{g}\right)$ центр
алгебры Ли $\frak{g}$, а через $b_{1}(\frak{g})$ --- первое число
Бетти, то есть размерность пространства
$H_{1}(\frak{g},\mathbb{R})$.

\begin{theorem}
Пусть $\frak{g}_{i}$ $\left( i=1,2\right) $ --- разрешимые алгебры Ли, пусть $%
b_{1}\left( \frak{g}_{i}\right) =\dim \left(
\frac{\frak{g}_{i}}{\left[ \frak{g}_{i},\frak{g}_{i}\right]
}\right) $. Тогда существует единственное центральное расширение
$\frak{e}$ для $\frak{g}_{1}\oplus \frak{g}_{2}$, удовлетворяющее
следующим условиям:

\begin{enumerate}
\item  $b_{1}\left( \frak{e}\right) =b_{1}\left( \frak{g}_{1}\right)
+b_{1}\left( \frak{g}_{2}\right) $

\item  $\dim D^{1}\left( \frak{e}\right) =\dim D^{1}\left( \frak{g}%
_{1}\right) +\dim D^{1}\left( \frak{g}_{2}\right) +b_{1}\left( \frak{g}%
_{1}\right) b_{1}\left( \frak{g}_{2}\right) $

\item  $\dim D^{i}\left( \frak{e}\right) =\dim D^{i}\left( \frak{g}%
_{1}\right) +\dim D^{i}\left( \frak{g}_{2}\right) ,\;i\geq 2$

\item  $\dim \,Z\left( \frak{e}\right) =\dim Z\left( \frak{g}_{1}\right)
+\dim Z\left( \frak{g}_{2}\right) +b_{1}\left( \frak{g}_{1}\right)
b_{1}\left( \frak{g}_{2}\right) .$
\end{enumerate}
\end{theorem}

\begin{proof}
Доказательство --- формально то же самое, как и в нильпотентном
случае \cite{C2,AC4}, при этом убывающая центральная
последовательность заменяется на производую последовательность.
Фактически расширение определяется посредством коцикла $\varphi \in H^{2}\left( \frak{g}_{1}\oplus \frak{g}_{2},\mathbb{R}^{b_{1}\left( \frak{g}%
_{1}\right) b_{2}\left( \frak{g}_{2}\right) }\right) $, так что
\begin{equation}
\varphi \left( \frak{g}_{1},D^{1}\frak{g}_{2}\right) =\varphi \left( D^{1}%
\frak{g}_{1},\frak{g}_{2}\right) =0.
\end{equation}
Условие (1) подразумевает, что
\begin{equation}
\dim _{\mathbb{R}}\left\{ \varphi \left( X,X^{\prime }\right) \;|\;X\in
\frak{g}_{1}-\left[ \frak{g}_{1},\frak{g}_{1}\right] ,\;X^{\prime }\in \frak{%
g}_{2}-\left[ \frak{g}_{2},\frak{g}_{2}\right] \right\} =b_{1}\left( \frak{g}%
_{1}\right) b_{1}\left( \frak{g}_{2}\right) .
\end{equation}
Без ограничения общности мы можем считать, что $\left\{
X_{1},..,X_{b_{1}\left( \frak{g}_{1}\right) }\right\} $ и $\left\{
X_{1}^{\prime },..,X_{b_{1}\left( \frak{g}_{2}\right) }^{\prime
}\right\} $ --- минимальные наборы образующих для $\frak{g}_{1}$ и
$\frak{g}_{2}$ соответственно.
Если $\left\{ e_{1},..,e_{b_{1}\left( \frak{g}_{1}\right) b_{1}\left( \frak{g%
}_{2}\right) }\right\} $ --- канонический базис в
$\mathbb{R}^{b_{1}\left( \frak{g}_{1}\right) b_{2}\left(
\frak{g}_{2}\right) }$, то мы можем выбрать такого представителя
$\varphi $, для которого имеет место
\begin{eqnarray}
\varphi \left( X_{i},X_{j}^{\prime }\right)  &=&e_{\left( i-1\right)
b_{1}\left( g_{2}\right) +j},\;1\leq i\leq b_{1}\left( \frak{g}_{1}\right)
,\;1\leq j\leq b_{1}\left( \frak{g}_{2}\right)  \\
\varphi \left( X,X^{\prime }\right)  &=&0\text{ при }X\in D^{1}\frak{g}_{1}%
\text{ или }X^{\prime }\in D^{1}\frak{g}_{2}
\end{eqnarray}
Соответствующее расширение удовлетворяет всем требованиям теоремы.
\end{proof}

Пусть $B_{i}$ --- базис алгебры $\frak{g}_{i}$ ($i=1,2$). Тогда
скобки Ли в алгебре $\frak{g}_{1}\underline{\times }\frak{g}_{2}$
представляют собой объединение скобок для алгебр $\frak{g}_{1}$ и
$\frak{g}_{2}$ и, кроме того, скобок вида $[X,X^{\prime}]$ для
всех образующих $X\in \frak{g}_{1}$ и $X^{\prime }\in
\frak{g}_{2}$. В силу уравнения Маурера-Картана, искомое
расширение может быть получено добавлением к  $b_{1}\left(
\frak{g}_{1}\right) b_{1}\left( \frak{g}_{2}\right) $ 2-формы вида
\begin{equation*}
d\eta _{ij}=\omega _{i}\wedge \omega _{j}^{\prime },
\end{equation*}
где $\omega _{i}$ (соотв. $\omega _{i}^{\prime }$) --- форма,
соответствующая образующей $X_{1}\in \frak{g}_{1}$ (соотв.,
$X_{2}\in \frak{g}_{2}$). Определим теперь произведение
образующими следующим образом.

\begin{definition}
Произведение образующими $\frak{g}_{1}\underline{\times
}\frak{g}_{2}$ представляет собой центральное расширение $\
\frak{g}_{1}\oplus \frak{g}_{2}$, определенное классом гомологий
$\varphi $.
\end{definition}

\begin{proposition}
Если $\frak{g}_{i}$ --- разрешимая алгебра Ли индекса $j_{i}$ при
$i=1,2$, то $\frak{g}_{1}\underline{\times }\frak{g}_{2}$ ---
неразложимая алгебра Ли индекса  $\max \left\{ j_{1},j_{2}\right\} $. В частности, $\frak{%
g}_{1}\underline{\times }\frak{g}_{2}=\frak{g}_{2}\underline{\times }\frak{g}%
_{1}.$
\end{proposition}

То есть, расширение $\frak{%
g}_{1}\underline{\times }\frak{g}_{2}=\frak{g}_{2}\underline{\times }\frak{g}%
_{1}$ является разрешимой алгеброй Ли.

\begin{definition}
Пусть $\frak{g}$ --- разрешимая алгебра Ли. Алгебра $\frak{g}$
называется произведением образующими, если существуют подалгебры
$\frak{g}_{1}$ и $\frak{g}_{2}$, такие что $\frak{g\simeq
g}_{1}\underline{\times }\frak{g}_{2}$.
\end{definition}

\begin{proposition}
Для любого  $d\geq 4$ существует разрешимая алгебра Ли $d$-го
порядка. Тогда существует пара $(\frak{g}_{1},\frak{g}_{2})$
(здесь $\frak{g}_{i}$ --- подалгебра), т.ч. $\frak{g\simeq
g}_{1}\underline{\times }\frak{g}_{2}$.
\end{proposition}

\begin{proof}
Для любого  $n\geq 3$ имеем разбиение
\begin{equation}
n=\left( n-2\right) +1+1.
\end{equation}
Пусть $\frak{r}_{n}$ --- алгебра Ли, в которой скобка имеет вид
\begin{equation}
\left[ X_{1},X_{j}\right] =X_{j},\;2\leq n-2
\end{equation}
в базисе $\left\{ X_{1},..,X_{n-2}\right\} $. Алгебра Ли
$\frak{r}_{n}$ --- разрешима и не нильпотентна. Далее, $\dim
\frak{r}_{n}/\left[ \frak{r}_{n},\frak{r}_{n}\right] =1$. Если
$L_{1}$ --- абелева алгебра первого порядка, то произведение
$\frak{r}_{1}\underline{\times }{L}_{1}$ представляет собой
разрешимую алгебру $n$-го порядка (за исключением случая $n=3$,
так как $\frak{r}_{3}=L_{1}$, следовательно
$\frak{r}_{3}\underline{\times }{L}_{1}$ изоморфна алгебре
Гейзенберга $\frak{h}_{1}$).
\end{proof}

Если одна из алгебр $\frak{g}_{1}$ или $\frak{g}_{2}$ разрешима,
то произведение этих алгебр также разрешима. Имеют место следующие
свойства:

\begin{proposition}
Пусть $\frak{g}$ --- произведение образующими. Тогда
\begin{enumerate}
\item  $\dim \left[ \frak{g},\frak{g}\right] \leq \dim \frak{g}-2$.

\item  $\dim Z\left( \frak{g}\right) \geq 1.$

\item  $\dim \frak{g}\geq 3.$

\item  Если алгебра $\frak{g}$ разрешима и ненильпотентна, то
существует разрешимая ненильпотентная подалгебра $\frak{g}_{1}$ и
$\frak{g}_{2}$ такая, что $\frak{g\simeq g}_{1}\underline{\times
}\frak{g}_{2}$.
\end{enumerate}
\end{proposition}

Принимая во внимание перечисленные выше структурные результаты о
произведениях образующими, мы можем получить полную классификацию
таких алгебр, опираясь на уже известную классификацию в низших
размерностях. Приведем классификаицию в размерности семь с
точностью до изоморфизма.

\begin{proposition}
Пусть $\frak{g\simeq }$ $\frak{g}_{1}\underline{\times }\frak{g}_{2}$ ---
произведение образующими, и пусть $\left( d_{i},m_{i}\right) =\left( \dim \frak{g}_{i},b_{1}\left( \frak{g}_{i}\right) \right) $
при $i=1,2$. Если $\dim \frak{g}%
=7$, то имеет место один из следующих случаев:

\begin{enumerate}
\item  $\left( d_{1},m_{1}\right) =\left( 5,1\right) $, $\left(
d_{2},m_{2}\right) =\left( 1,1\right) $.

\item  $\left( d_{1},m_{1}\right) =\left( 4,1\right) $, $\left(
d_{2},m_{2}\right) =\left( 2,1\right) $.

\item  $\left( d_{1},m_{1}\right) =\left( 3,1\right) $, $\left(
d_{2},m_{2}\right) =\left( 3,1\right) $.

\item  $\left( d_{1},m_{1}\right) =\left( 4,2\right) $, $\left(
d_{2},m_{2}\right) =\left( 1,1\right) $.

\item  $\left( d_{1},m_{1}\right) =\left( 3,2\right) $, $\left(
d_{2},m_{2}\right) =\left( 2,1\right) $.

\item  $\left( d_{1},m_{1}\right) =\left( 3,1\right) $, $\left(
d_{2},m_{2}\right) =\left( 2,2\right) $

\item  $\left( d_{1},m_{1}\right) =\left( 3,3\right) $, $\left(
d_{2},m_{2}\right) =\left( 1,1\right) $.
\end{enumerate}

В частности, если $b_{1}\left( \frak{g}\right) =4$, то алгебра
$\frak{g}$ является нильпотентной.
\end{proposition}

\begin{proof}
Так как $\dim \left( \frak{g}_{1}\underline{\times
}\frak{g}_{2}\right) =d_{1}+d_{2}+m_{1}m_{2}$, дальнейшая
классификация сводится к разбиению числа семь в сумму трех
натруальных чисел. Учитывая $b_{1}\left(
\frak{g}_{1}\underline{\times }\frak{g}_{2}\right) =b_{1}\left(
\frak{g}_{1}\right) +b_{1}\left( \frak{g}_{2}\right) =m_{1}+m_{2}$
и $\dim \frak{g}=7$, мы имеем $2\leq b_{1}\left(
\frak{g}_{1}\underline{\times }\frak{g}_{2}\right) \leq 4$. Если
$b_{1}\left( \frak{g}_{1}\underline{\times }\frak{g}_{2}\right)
=2$, to $m_{1}=m_{2}=1$, и мы приходим к случаям 1.-3. Если $
b_{1}\left( \frak{g}_{1}\underline{\times }\frak{g}_{2}\right)
=3$, мы имеем
\begin{equation*}
7=d_{1}+d_{2}+2.
\end{equation*}
Если $d_{1}=4$, то $d_{2}=1$, следовательно $m_{2}=1$, $m_{2}=2$.
Мы приходим к случаю 4. Если  $d_{1}=3$ и $d_{2}=2$, то либо
$m_{1}=1$ и $m_{2}=2$, либо $m_{1}=2$ и $m_{2}=1$. Это приводит
нас к случаям 5 и 6. Если $b_{1}\left(
\frak{g}_{1}\underline{\times }\frak{g}_{2}\right) =4$, то
$\frak{g}_{1}$ и $\frak{g}_{2}$ --- абелевы алгебры, и
произведение $\frak{g}_{1}\underline{\times }\frak{g}_{2}$
нильпотентно.
\end{proof}

Для порядков $n\geq 8$ мы можем решать эту задачу, но так как
разрешимые алгберы Ли размерностей $n\geq7$ не классифицированы
(последняя известная размерность --- шесть \cite{Mu4,Tu}),
результат не был бы окончательным.

\section{Применение к инвариантам}

Инварианты коприсоединенного представления играют важную роль в
теории представлений, а также во многих геометрических и
физических приложениях (квантовые числа и ограничения
представлений для моделей в ядерной и атомной физике).\newline

Одним из стандартных методов для определения этих инвариантов
является теория дифференциальных уравнений; она позволяет получать
не только полиномиальные инварианты (связанные с операторами
Казимира), но также и иррациональные функции, которые находят
применение в теории интегрируемых гамильтоновых систем
\cite{Pe,Ber,Tro}. Пусть $\left\{C_{ij}^{k}\right\}$ ---
структурный тензор алгебры Ли $\frak{g}$ в базисе
$\left\{X_{1},..,X_{n}\right\}  $. Пусть
\begin{equation}
\widehat{X}_{i}=-C_{ij}^{k}x_{k}\partial_{x_{j}}
\end{equation}
--- представление алгебры Ли, где $\left[  X_{i},X_{j}\right]  =C_{ij}^{k}X_{k}$ $\left(  1\leq i<j\leq
n,\;1\leq k\leq n\right)  $. Аналитическая функция $F\left(
X_{1},..,X_{n}\right) $ на $\frak{g}^{*}$ является инвариантом
коприсоединенного представления алгебры Ли $\frak{g}$ тогда и
только тогда, когда
\begin{equation}
\widehat{X}_{i}F\left(  x_{1},..,x_{n}\right)  =-C_{ij}^{k}x_{k}%
\partial_{x_{j}}F\left(  x_{1},..,x_{n}\right)  =0,\;1\leq i\leq n. \label{sys}
\end{equation}
Общее число  $\mathcal{N}\left(  \frak{g}\right)  $ решений (как
полиномиальных, так и неполиномиальных) системы (\ref{sys})
задается формулой
\begin{equation}
\mathcal{N}\left(  \frak{g}\right)  =\dim\,\frak{g}-{\rm rank} \,\left(C_{ij}^{k}x_{k}\right),
\end{equation}
где матрица  $A(\frak{g})$ выражается по формуле:
\begin{equation}
A(\frak{g})=\left(C_{ij}^{k}x_{k}\right),
\end{equation}
в базисе $\left\{X_{1},..,X_{n}\right\}  $. Эту формулу можно
также интерпретировать с использованием уравнения Маурера-Картана
следующим образом. Пусть $\mathcal{L}(\frak{g})=\mathbb{R}\left\{
d\omega_{i}\right\} _{1\leq i\leq \dim\frak{g}}$ ---
подпространство пространства $\bigwedge^{2}\frak{g}^{\ast}$,
порожденное 2-формами $d\omega_{i}$. Так как
$\dim\mathcal{L}(\frak{g})=\dim\left( \frak{g}\right)  $ тогда и
только тогда, когда $d\omega_{i}\neq0$ при всех $i$, то есть
$\frak{g}=[\frak{g},\frak{g}]$, если $\omega=a^{i}d\omega_{i}$ ---
элемент общего положения в  $\mathcal{L}(\frak{g})$, тогда
существует $j_{0}\left( \omega\right) \in\mathbb{N}$, т.ч.
\begin{equation}
\bigwedge^{j_{0}\left(  \omega\right)  }\omega\neq0,\quad \bigwedge
^{j_{0}\left(  \omega\right)  +1}\omega\equiv0.
\end{equation}
Тогда $r\left(  \omega\right)  =2j_{0}\left(  \omega\right)$, где
$r\left(  \omega\right)$ --- ранг 2-формы $\omega$ \cite{C43}.
Положим
\begin{equation}
j_{0}\left(  \frak{g}\right)  =\max\left\{  j_{0}\left(  \omega\right)
\;|\;\omega\in\mathcal{L}(\frak{g})\right\}.
\end{equation}
Очевидно, что  $j_{0}\left(  \frak{g}\right)  $ --- числовой
инвариант алгебры Ли $\frak{g}$. Короткое вычисление показывает,
что
\begin{equation}
\mathcal{N}\left(  \frak{g}\right)  =\dim\frak{g}-2j_{0}\left(  \frak{g}\right) \label{Ca}.
\end{equation}
Формула (\ref{Ca}) важна для оценки числа инвариантов полупрямых
произведений в том случае, когда общий метод уже не работает.

Для определения числа независимых инвариантов применим формулу
(\ref{Ca}), учитывая уравнения Маурера-Картана. Пусть
$\frak{g}_{1}$ и $\frak{g}_{2}$ --- разрешимые алгебры Ли и пусть
$\left\{ C_{ij}^{k}\right\}$ (соответстсвенно $\left\{
\overline{C}_{lm}^{n}\right\} $) --- структурный тензор алгебры
$\frak{g}_{i}$ ($i=1,2$) в базисе $\left\{
X_{1},..,X_{d_{1}}\right\} $, соотв. $\left\{ X_{1}^{\prime
},..,X_{d_{2}}^{\prime }\right\} $. Уравнения Маурера-Картана
имеют вид:
\begin{eqnarray}
d\omega _{i} &=&0,\;1<i<m_{1} \label{MC1}  \\
d\omega _{k} &=&-C_{ij}^{k}\omega _{i}\wedge \omega _{j},\;k\geq m_{1},
\end{eqnarray}
для $\frak{g}_{1}$, при этом
\begin{eqnarray}
d\omega _{i}^{\prime } &=&0,\;1<i<m_{2} \\
d\omega _{k}^{\prime } &=&-\overline{C}_{ij}^{k}\omega _{i}^{\prime }\wedge
\omega _{j}^{\prime },\;k\geq m_{2}. \label{MC4}
\end{eqnarray}
для $\frak{g}_{2}$. Тогда уравнения для произведений
$\frak{g}_{1}\underline{\times }\frak{g}_{2}$  имеют вид (15)-(18)
i
\begin{equation}
d\eta_{ij}=\omega_{i}\wedge\omega_{j}^{\prime},\quad 1\leq i\leq m_{1}, 1\leq j\leq m_{2}. \label{MC5}
\end{equation}
Для всех $X_{i}$ (соответственно для $X_{i}^{\prime }$)
\begin{equation}
X_{i}\lrcorner \,\omega _{j}=\delta _{ij}\;(\text{соотв.
}X_{i}^{\prime }\lrcorner \,\omega _{j}^{\prime }),
\end{equation}
где $\lrcorner$ --- внутренне произведение. Число
$\mathcal{N}\left(  \frak{g}\right)$ можно получить из
максимального ранга элементов общего положения, задаваемого
формулами (15)-(19). Приводимый ниже результат задает метод
нахождения этого числа.

\begin{proposition}
Пусть  $\omega _{0}\in \mathcal{L}\left( \frak{g}_{1}\right) $, и
пусть $\omega
_{0}^{\prime }\in \mathcal{L}\left( \frak{g}_{2}\right) $ --- две формы, т.ч.
$j_{0}\left( \omega _{0}\right) =j_{0}\left( \frak{g}_{1}\right) $ и $%
j_{0}\left( \omega _{0}^{\prime }\right) =j_{0}\left(
\frak{g}_{2}\right) $. Пусть $F\subset \left\{ 1,..,m_{1}\right\}
$ и $F^{\prime }\subset\left\{ 1,..,m_{2}\right\} $ ---
подмножества, т.ч.
\begin{enumerate}
\item  $\forall \alpha \in F$ $\Rightarrow X_{\alpha }\lrcorner \,\omega
_{0}\equiv 0.$

\item  $\forall \beta \in F^{\prime }\Rightarrow X_{\beta }^{\prime
}\lrcorner \,\omega _{0}^{\prime }\equiv 0$.
\end{enumerate}
Тогда
\begin{equation}
j_{0}\left( \frak{g}_{1}\underline{\times }\frak{g}_{2}\right) =j_{0}\left(
\frak{g}_{1}\right) +j_{0}\left( \frak{g}_{2}\right) +j_{0}\left( \omega
\right) ,
\end{equation}
где
\begin{equation}
\omega =\sum_{\alpha \in F,\beta \in F^{\prime }}\omega _{\alpha }\wedge
\omega _{\beta }.
\end{equation}
В частности,
\begin{equation}
m_{1}m_{2}\leq \mathcal{N}\left( \frak{g}_{1}\underline{\times }\frak{g}%
_{2}\right) \leq \mathcal{N}\left( \frak{g}_{1}\right) +\mathcal{N}\left(
\frak{g}_{2}\right) +\left( m_{1}m_{2}-2j_{0}\left( \omega \right) \right) .
\end{equation}
\end{proposition}

\begin{proof}
Уравнния Маурера-Картана  $\frak{g}_{1}\underline{\times
}\frak{g}_{2}$ получаются из уравнений (15)-(18). Число
независимых инвариантов выписывается по формуле
\begin{equation}
\mathcal{N}\left( \frak{g}_{1}\underline{\times }\frak{g}_{2}\right) =\dim
\left( \frak{g}_{1}\underline{\times }\frak{g}_{2}\right) -2j_{0}\left(
\theta \right) ,
\end{equation}
где $\theta \in \mathcal{L}\left( \frak{g}_{1}\underline{\times }\frak{g}%
_{2}\right) $ --- общая форма максимального ранга. Пусть $\alpha \in F$. Так как $%
\omega _{0}\in \mathcal{L}\left( \frak{g}_{1}\right) $, уравнение $%
X_{\alpha }\,\lrcorner \,\omega _{0}=0$ подразумевает, что $\left[ X,X_{\alpha }%
\right] =\psi ^{k}X_{k}\neq 0$ для всех $X\in \frak{g}_{1}$ (то
есть, $d\omega _{k}$ не лежит в $\omega_{0}$). Легко видеть, что
\begin{equation}
\omega _{0}+\omega _{0}^{\prime }+\sum_{\alpha \in F,\beta \in F^{\prime
}}\omega _{\alpha }\wedge \omega _{\beta }
\end{equation}
имеет максималный ранг, так как $j_{0}\left( \omega _{0}\right) =j_{0}\left( \frak{g}%
_{1}\right) ,\;j_{0}\left( \omega _{0}^{\prime }\right) =j_{0}\left( \frak{g}%
_{2}\right) $ и $\omega \notin \mathcal{L}\left(
\frak{g}_{1}\right) + \mathcal{L}\left( \frak{g}_{2}\right) $.
Тогда
\begin{eqnarray}
\mathcal{N}\left( \frak{g}_{1}\underline{\times }\frak{g}_{2}\right)
&=&\dim \frak{g}_{1}+\dim \frak{g}_{2}+m_{1}m_{2}-2\left( j_{0}\left( \omega
_{0}\right) +j_{0}\left( \omega _{0}^{\prime }\right) +j_{0}\left( \omega
\right) \right) \nonumber \\
&=&\mathcal{N}\left( \frak{g}_{1}\right) +\mathcal{N}\left( \frak{g}%
_{2}\right) +\left( m_{1}m_{2}-2j_{0}\left( \omega \right) \right) .
\end{eqnarray}
$m_{1}m_{2}\leq \mathcal{N}\left( \frak{g}_{1}\underline{\times
}\frak{g }_{2}\right) $ следует из того факта, что добавленные
элементы являются центральными.
\end{proof}

\begin{corollary}
Если $\mathcal{N}\left( \frak{g}_{1}\right) =\mathcal{N}\left( \frak{g%
}_{2}\right) =0$, то $\mathcal{N}\left(
\frak{g}_{1}\underline{\times }\frak{g}_{2}\right) =m_{1}m_{2}$.
\end{corollary}

\begin{proof}
В этом случае  $\mathcal{N}\left( \frak{g}_{1}\right) =\mathcal{N}\left( \frak{g}%
_{2}\right) =0$, существуют $2$-формы $\omega _{0}\in \mathcal{L}\left( \frak{%
g}_{1}\right) $ и $\omega _{0}^{\prime }\in \mathcal{L}\left( \frak{g}%
_{2}\right) $, т.ч. $\dim \frak{g}_{1}=2j_{0}\left( \omega
_{0}\right) $ и $\dim \frak{g}_{2}=2j_{0}\left( \omega
_{0}^{\prime }\right) $. То есть, форма
$\theta=\omega_{0}+\omega_{0}^{\prime}$ имеет ранг $\dim
\frak{g}_{1}+\dim \frak{g}_{2}$. Алгебра $
\frak{g}_{1}\underline{\times }\frak{g}_{2}$ содержит в точности
$m_{1}m_{2}$ центральных элементов. Таким образом, форма $\theta$
имеет максимальный ранг.
\end{proof}

\begin{remark}
Обратное, вообще говоря, неверно. Пусть $\frak{g}$ --- алгебра Ли,
задаваемая уравнениями
\begin{eqnarray*}
d\omega _{1} &=&\omega _{2}\wedge \omega _{3}+\omega _{1}\wedge \omega _{4}
\\
d\omega _{2} &=&\omega _{2}\wedge \omega _{4}-\omega _{2}\wedge \omega _{5}
\\
d\omega _{3} &=&\omega _{3}\wedge \omega _{5} \\
d\omega _{4} &=&d\omega _{5}=0
\end{eqnarray*}
Форма $\omega _{0}=d\omega _{1}$ удовлетворяет уравнению
$j_{0}\left( \omega
_{0}\right) =j_{0}\left( \frak{g}\right) =2$, то есть $\mathcal{N}\left( \frak{g%
}\right) =1$. Тогда $\frak{g}\underline{\times }\frak{g}$ ---
разрешимая алгебра Ли 14-го порядка (так как $b_{1}\left(
\frak{g}\right) =2$), при этом уравнения Маурера-Картана имеют вид
\begin{equation*}
\begin{tabular}{lll}
$d\omega _{1}=\omega _{2}\wedge \omega _{3}+\omega _{1}\wedge \omega _{4};$
& $d\omega _{1}^{\prime }=\omega _{2}^{\prime }\wedge \omega _{3}^{\prime
}+\omega _{1}^{\prime }\wedge \omega _{4}^{\prime };$ & $d\eta _{1}=\omega
_{4}\wedge \omega _{4}^{\prime }$ \\
$d\omega _{2}=\omega _{2}\wedge \omega _{4}-\omega _{2}\wedge \omega _{5};$
& $d\omega _{2}^{\prime }=\omega _{2}^{\prime }\wedge \omega _{4}^{\prime
}-\omega _{2}^{\prime }\wedge \omega _{5}^{\prime };$ & $d\eta _{2}=\omega
_{4}\wedge \omega _{5}^{\prime }$ \\
$d\omega _{3}=\omega _{3}\wedge \omega _{5};$ & $d\omega _{3}^{\prime
}=\omega _{3}^{\prime }\wedge \omega _{5}^{\prime };$ & $d\eta _{3}=\omega
_{5}\wedge \omega _{4}^{\prime }$ \\
$d\omega _{4}=d\omega _{5}=0;$ & $d\omega _{4}^{\prime }=d\omega
_{5}^{\prime }=0;$ & $d\eta _{4}=\omega _{5}\wedge \omega _{5}^{\prime }$%
\end{tabular}
\end{equation*}
в базисе $\left\{ \omega _{1},..,\omega _{5},\omega _{1}^{\prime
},..,\omega _{5}^{\prime },\eta _{1},..,\eta _{4}\right\} $. Пусть
$\theta =\alpha ^{i}d\omega _{i}+\beta ^{j}d\omega _{j}^{\prime
}+\gamma
^{k}d\eta _{k}\in \mathcal{L}\left( \frak{g}\underline{\times }\frak{g}%
\right) $ --- элемент общего положения.Тогда $\bigwedge^{6}\theta
\equiv 0$ и
\begin{equation*}
\bigwedge^{5}\theta =\left( \alpha ^{1}\right) ^{2}\left( \beta ^{1}\right)
^{2}\gamma ^{4}\omega _{1}\wedge \omega _{2}\wedge \omega _{3}\wedge \omega
_{4}\wedge \omega _{5}\wedge \omega _{1}^{\prime }\wedge \omega _{2}^{\prime
}\wedge \omega _{3}^{\prime }\wedge \omega _{4}^{\prime }\wedge \omega
_{5}^{\prime }\neq 0
\end{equation*}
тогда и только тогда, кода $\alpha ^{1}\beta ^{1}\gamma ^{4}\neq
0$. Форма $\theta ^{\prime }=d\omega _{1}+d\omega _{1}^{\prime
}+d\eta _{4}$ имеет максимальный ранг, при этом $\mathcal{N}\left(
\frak{g}\underline{\times }\frak{g}\right) =4$. В частности, мы
получаем, что если $F$ и $F^{\prime}$ не пусты, то этом $\omega
=d\eta _{4}$ удовлетворяет условию $j_{0}\left( \omega \right)
=1$, тогда $2\mathcal{N}\left( \frak{g}\right) -2j_{0}\left(
\omega \right) =0$. Кроме того, алгебра $\frak{g}\underline{\times
}\frak{g}$ имеет четыре оператора Казимира, в то же время
$\frak{g}$ не имеет полиномиальных инвариантов (но лишь
рациональные инварианты).

\end{remark}

\section{Приложения к структурам произведений}

Приведем теперь некоторые приложения к
дифференциально-геометрическим структурам на группах Ли. Эти
конструкции позволяют получать некоторые новые алгебры Ли в
произвольных размерностях, которые образуют новые классы,
пополняющие известные случаи.

\begin{definition}
Пусть $E:\frak{g\rightarrow g}$, $E\neq \pm id$ --- инволютивный
автоморфизм алгебры $\frak{g}$. В этом случае автоморфизм $E$
задает структуру почти произведения. Говорят, что $E$ задает
структуру произведения, если имеет место интегрируемость, т.е.
если для любых $X,Y\in \frak{g}$ имеют место следующие тождества.
\begin{equation*}
E\left[ X,Y\right] =\left[ E\left( X\right) ,Y\right] +\left[ X,E\left(
Y\right) \right] -E\left[ E\left( X\right) ,E\left( Y\right) \right] .
\end{equation*}
\end{definition}

Мы имеем разложение
\begin{equation}
\frak{g=g}_{+}\oplus \frak{g}_{-}, \label{DD}
\end{equation}
где $\frak{g}_{+}$ --- собственное подпространство,
соответствующее собственному значению $\lambda=1$, а
$\frak{g}_{-}$ --- собственное пространство, соответствующее
собственному значению $\lambda=-1$. Если условие интегрируемости
выполняется, то подпространства $\frak{g}_{+}$ и $\frak{g}_{-}$
являются подалгебрами. Такие структуры были изучены для различных
типов алгебр Ли и классифицированы для алгебр Ли размерностей
$n\leq 4$ \cite{Ov}. В настоящем разделе мы показываем, что
произведения образующими дают новые нетривиальным примеры алгебр
Ли, обладающие такой структурой; эти алгебры определяются в
терминах своих подалгебр.

\begin{proposition}
Пусть $\frak{g}_{i}$ $\left( i=1,2\right) $ --- разрешимые алгебры
Ли, и пусть $E_{i}$ $(i=1,2)$ --- структура произведения для
$i=1,2$. Пусть $E_{1}\left( X_{j}\right) =\varepsilon _{j}X_{j}$
при $1\leq j\leq n\,$, $\varepsilon _{j}=\pm 1$ в базисе $\left\{
X_{1},..,X_{n}\right\} $ (соотв. $E_{2}\left(
X_{j}^{\prime}\right) =\varepsilon _{j}X_{j}^{\prime}$ при $1\leq
j\leq m\,$, $\varepsilon _{j}=\pm 1$ в базисе $\left\{
X_{1}^{\prime},..,X_{n}^{\prime}\right\} $). Тогда
$\frak{g}_{1}\underline{%
\times }\frak{g}_{2}$ является структурным произведением.
\end{proposition}

\begin{proof}
Без ограничения общности будем считать, что $X_{i}\left( 1\leq
i\leq b_{1}\left( \frak{g}_{1}\right) \leq n\right) $, т.ч.
$X_{i}\notin \left[ \frak{g}_{1},\frak{g}_{1}\right] $ (по
аналогии со случаем $\frak{g}_{2}$). Пусть $\left[
X_{i},X_{j}^{\prime }\right] =Z_{ij}$
при  $1\leq i\leq b_{1}(%
\frak{g}_{1}),1\leq j\leq b_{1}\left( \frak{g}_{2}\right) $ и $E:\frak{%
g}_{1}\underline{\times }\frak{g}_{2}\rightarrow \frak{g}_{1}\underline{%
\times }\frak{g}_{2}$ by $E|_{\frak{g}_{1}}=E_{1}$ --- векторное
отображение, т.ч. $E|_{\frak{g}_{1}}=E_{1}$ и пусть
$E|_{\frak{g}_{2}}=E_{2}$. Кроме того, действие $E$ на новых
элементах $Z_{ij}$ --- диагонально. Так как  $E_{1}$ и $E_{2}$
интегрируемы, при этом $\left[
\frak{g}_{1},D^{1}\frak{g}_{2}\right] =\left[ D^{1}\frak{g}_{1},\frak{g}_{2}%
\right] =0$, нам достаточно проверить условие интегрируемости лишь
на парах $\left( X,X^{\prime }\right) \in \frak{g}_{1}\frak{\oplus
g}_{2}$ c $\left[ X,X^{\prime }\right] \neq 0$: Пусть $X=\alpha
^{i}X_{i}+\beta ^{j}Y_{j}
$ $\left( i\leq n\right) $ и $X^{\prime }=\gamma ^{k}X_{k}^{\prime}+\delta ^{l}Y_{l}^{\prime }$, где $Y_{i}\in \left[ \frak{g}_{1},\frak{g}%
_{1}\right] $ соответственно  $Y_{i}^{\prime }\in \left[ \frak{g}_{2},\frak{g}%
_{2}\right] $. Тогда
\begin{eqnarray}
&&E\left[ X,X^{\prime }\right] -\left[ E\left( X\right) ,X^{\prime }\right] -%
\left[ X,E\left( X^{\prime }\right) \right] +E\left[ E\left( X\right)
,E\left( X^{\prime }\right) \right] \nonumber \\
&=&\alpha ^{i}\gamma ^{k}E\left( Z_{ik}\right) -\alpha ^{i}\gamma ^{k}\left(
\varepsilon _{i}+\varepsilon _{k}\right) Z_{ik}+\alpha ^{i}\gamma
^{k}\varepsilon _{i}\varepsilon _{k}E\left( Z_{ik}\right) \nonumber \\
&=&\alpha ^{i}\gamma ^{k}\left\{ \left( 1+\varepsilon _{i}\varepsilon
_{k}\right) E\left( Z_{ik}\right) -\left( \varepsilon _{i}+\varepsilon
_{k}\right) Z_{ik}\right\}. \label{AT}
\end{eqnarray}

Чтобы уравнение (28) обратилось в ноль, для всех пар $(i,k)$ нам
нужно $E\left( Z_{ik}\right) =\lambda _{ik}Z_{ik}$ s $\lambda
_{ik}=\pm 1$. Имеют место два различных случая:

\begin{enumerate}
\item  При $\left( i_{0},k_{0}\right) $ мы имеем $\varepsilon
_{i_{0}}+\varepsilon _{k_{0}}=0$. Тогда $1+\varepsilon
_{i_{0}}\varepsilon _{k_{0}}=0$, и мы можем положить  $\lambda
_{i_{0}k_{0}}$ равным $1$ или $-1$.

\item  При $\left( i_{0},k_{0}\right) $ имеем $\varepsilon
_{i_{0}}+\varepsilon _{k_{0}}\neq 0$. Тогда $sgn\left( \varepsilon
_{i_{0}}\right) =sgn\left( \varepsilon _{k_{0}}\right) $; таким
образом,
 $2E\left( Z_{i_{0}k_{0}}\right) =\pm 2Z_{i_{0}k_{0}}$, поэтому $%
\lambda _{i_{0}k_{0}}=\pm 1$; знак выбирается в зависимости от
$\varepsilon _{i_{0}}+\varepsilon _{k_{0}}$.
\end{enumerate}

Действие $E$ на $Z_{ij}$ полностью определено в случае $\left(
\varepsilon _{i}+\varepsilon _{k}\right) \neq 0$. В остальных
случаях мы можем выбрать знак для $\lambda_{ij}$. Для разложения
(\ref{DD}) мы имеем
\begin{eqnarray}
\frak{g}_{+}^{E} &=&\frak{g}_{+}^{E_{1}}\oplus \frak{g}_{+}^{E_{2}}\oplus
\left\langle Z_{ij}\;|\;\lambda _{ij}=1\right\rangle  \\
\frak{g}_{-}^{E} &=&\frak{g}_{-}^{E_{1}}\oplus \frak{g}_{-}^{E_{2}}\oplus
\left\langle Z_{ij}\;|\;\lambda _{ij}=-1\right\rangle
\end{eqnarray}
Так как $E_{1}$ и $E_{2}$ интегрируемы, то $\frak{g}%
_{+}^{E}$ и $\frak{g}_{-}^{E}$ являются подалгебрами в $\frak{g}_{1}\underline{%
\times }\frak{g}_{2}$.
\end{proof}

\begin{example}
Рассмотрим теперь следующие алгебры Ли $\frak{r}_{4}$ и
$\frak{r}_{4,0}$:
\begin{equation}
\frak{r}_{4}:\;\left[ X_{0},X_{1}\right] =X_{1},\;\left[ X_{0}X_{2}\right]
=X_{1}+X_{2},\;\left[ X_{0},X_{3}\right] =X_{2}+X_{3}
\end{equation}
в базисе $\left\{ X_{0},X_{1},X_{2},X_{3}\right\} $ и
\begin{equation}
\frak{r}_{4,0}:\;\left[ X_{4},X_{5}\right] =X_{6},\;\left[ X_{4},X_{7}\right]
=X_{6}
\end{equation}
в базисе $\left\{ X_{4},X_{5},X_{6},X_{7}\right\} $. Пусть
произведение $\frak{r}_{4}\underline{\times }\frak{r}_{4,0}$. В
базисе $\left\{ X_{0},..,X_{9}\right\} $ алгебра задается скобками
\begin{equation}
\left[ X_{0},X_{4}\right] =X_{8},\;\left[ X_{0},X_{7}\right] =X_{9}.
\end{equation}
и коммутационными соотношениями (31)-(32). Алгебры $\frak{r}_{4}$
и $\frak{r}_{4,0}$ имеют структуру $E_{1}$, соотв. $E_{2}$, где

\begin{eqnarray}
E_{1}\left( X_{i}\right)  &=&X_{i},\;i=0,1;\;E_{1}\left( X_{i}\right)
=-X_{i},\;i=2,3 \\
E_{2}\left( X_{i}\right)  &=&X_{i},\;i=4,5;\;E_{2}\left( X_{i}\right)
=-X_{i},\;i=6,7.
\end{eqnarray}
Теперь определим $E:\frak{r}_{4}\underline{\times }\frak{r}_{4,0}\rightarrow \frak{%
r}_{4}\underline{\times }\frak{r}_{4,0}$ таким образом, что $E|_{\frak{r}_{4}}=E_{1},\;E|_{%
\frak{r}_{4,0}}=E_{2}$. Из предложения 6 следует, что
\begin{equation}
E\left( X_{8}\right) =X_{8},\;E\left( X_{9}\right) =\lambda X_{9},
\end{equation}
где  $\lambda =\pm 1$. Если $ \lambda =1$, то мы получим структуру
на $\frak{r}_{4}\underline{\times }\frak{r}_{4,0}$, т.ч. $\dim
\left( \frak{r}_{4}\underline{\times }\frak{r}_{4,0}\right)
_{+}=6,\;\dim
\left( \frak{r}_{4}\underline{\times }\frak{r}_{4,0}\right) _{-}=4$, в то
время как выбор $\lambda =-1$ задает структуру с $\dim \left( \frak{r%
}_{4}\underline{\times }\frak{r}_{4,0}\right) _{+}=\dim \left( \frak{r}_{4}%
\underline{\times }\frak{r}_{4,0}\right) _{-}=5$.
\end{example}

\begin{definition}
Пусть  $\frak{g}$ --- алгебра Ли со структурой произведения $E$.
Если $\dim \frak{g}_{+}=\dim \frak{g}_{-}$, то $E$ называют
паракомплексной структурой на $\frak{g}$.
\end{definition}

\begin{corollary}
Пусть $\frak{g}$ --- разрешимая алгебра Ли, для которой $\left\{
X_{1},..,X_{4m},..,X_{2n}\right\} $, так что $b_{1}\left(
\frak{g}\right) =4m$, и пусть $E$ --- паракомплексная структура.
Если
\begin{enumerate}
\item  $X_{i}\notin D^{1}\frak{g}$ при $1\leq i\leq 4m,$

\item  $X_{i}\in \frak{g}_{+}$  при $1\leq i\leq 2m$ and $X_{i}\in \frak{g}%
_{-}$  for $2m+1\leq i\leq 4m.$

\item  $b_{1}\left( \frak{g}_{+}\right) =b_{1}\left( \frak{g}_{-}\right) =2m.
$
\end{enumerate}
то $\frak{g}_{+}\underline{\times }\frak{g}_{-}$ допускает
паракомплексную структуру.
\end{corollary}

\begin{proof}
Так как  $E$ --- паракомплексная структура, то размерность $\dim \frak{g}_{+}\underline{\times }\frak{g}_{-}$
четна. Предложение 6 задает структуру произведения $E$ на $\frak{g}_{+}%
\underline{\times }\frak{g}_{-}$. Ограничение на  $\frak{g}_{+}$ и
$\frak{g}_{-}$ является тождественным(соотв., антитождественным)
отображением. Условием 2 (см. доказательство предложения 6) знак
$E\left( \left[ X_{i},X_{j}\right] \right) $ $\left( 1\leq i\leq
2m,2m+1\leq j\leq 4m\right) $ определен быть не может. Из $4m^{2}$
элементов вида $[X_{i},X_{j}]$ сопоставим $2m^{2}$ элементам знак
плюс, а остальным  $2m^{2}$ элементам --- знак минус. Поэтому
$\dim \left( \frak{g}_{+}\underline{\times }\frak{g}_{-.}\right)
_{+}=\dim \left( \frak{g}_{+}\underline{\times
}\frak{g}_{-.}\right) _{-}=n+2m^{2}$.
\end{proof}

Предыдущий пример показывает, что условие является достаточным, но
не необходимым для существования комплексных структур на
произведениях.

Автор выражает благодарность А.В.Чернавскому за приглашение на
конференцию и В.О.Мантурову за ценные обсуждения и помощь.

}

\end{document}